\documentclass{amsart}%
\usepackage{amsfonts}
\usepackage{amsmath}
\usepackage{amssymb}
\usepackage{graphicx}%
\newtheorem{theorem}{Theorem}
\theoremstyle{plain}

\newtheorem{definition}{Definition}

\newtheorem{lemma}{Lemma}

\numberwithin{equation}{section}
\begin{document}
\title[Fractal transformations]{Transformations between attractors of hyperbolic iterated function systems}
\author{Michael F. Barnsley}
\address{Department of Mathematics\\
Australian National University\\
Canberra, ACT, Australia}
\email{michael.barnsley@maths.anu.edu.au, mbarnsley@aol.com}
\urladdr{http://www.superfractals.com}
\thanks{The author thanks Louisa Barnsley for help with the illustrations.}
\thanks{This paper is in final form and no version of it will be submitted for
publication elsewhere.}

\begin{abstract}
This paper is in the form of an essay. It defines fractal tops and code space
structures associated with set-attractors of hyperbolic iterated function
systems (IFSs). The fractal top of an IFS is associated with a certain shift
invariant subspace of code space, whence the entropy of the IFS, and of its
set-attractor, may be defined. Given any ordered pair of hyperbolic IFSs, each
with the same number of maps, there is a natural transformation, constructed
with the aid of fractal tops, whose domain is the attractor $A_{1}$of the
first IFS and whose range is contained in the attractor $A_{2}$ of the second
IFS. This transformation is continuous when the code space structure of the
IFS is "contained in" the code space structure of the second IFS, and is a
homeomorphism between $A_{1}$ and $A_{2}$ when the code space structures are
the same. Conversely, if two IFS are homeomorphic then they possess the same
code space structure. Hence we obtain that two IFS attractors are homeomorphic
then they have the same entropy. Several examples of fractal transformations
and fractal homeomphisms are given.

\end{abstract}
\date{March 13$^{th}$ 2007.}
\maketitle

\section{Introduction}

In this essay we introduce fractal transformations. The main examples are
fascinating mappings between diverse subsets of $\mathbb{R}^{2}$; they can be
readily illustrated by using the chaos game. Fractal transformations can be
quicky grasped because they rely on basic notions in topology, probability,
dynamical systems, and geometry. They may be applied to computer graphics to
produce digital content with new look-and-feel \cite{barnhutch}; they may also
be relevant to image compression and biological modelling.

\section{\label{hypersec}Hyperbolic IFS}

\begin{definition}
Let $\mathbb{(X},d_{\mathbb{X}})$ be a complete metric space. Let
$\{f_{1},f_{2},...,f_{N}\}$ be a finite sequence of strictly contractive
transformations, $f_{n}:\mathbb{X\rightarrow X}$, for $n=1,2,...,N$. Then
\[
\mathcal{F}:=\{\mathbb{X};f_{1},f_{2},...,f_{N}\}
\]
is called a hyperbolic iterated function system or hyperbolic IFS.
\end{definition}

A transformation $f_{n}:\mathbb{X\rightarrow X}$ is strictly contractive iff
there exists a number $l_{n}\in\lbrack0,1)$ such that $d(f_{n}(x),f_{n}%
(y))\leq l_{n}d(x,y)$ for all $x,y\in\mathbb{X}$. The number $l_{n}$ is called
a contractivity factor for $f_{n}$ and the number%
\[
l=\max\{l_{1},l_{2},...,l_{N}\}
\]
is called a contractivity factor for $\mathcal{F}$.

Let $\Omega$ denote the set of all infinite sequences of symbols $\{\sigma
_{k}\}_{k=1}^{\infty}$ belonging to the alphabet $\{1,...,N\}$. We write
$\sigma=\sigma_{1}\sigma_{2}\sigma_{3}...\in\Omega$ to denote a typical
element of $\Omega$, and we write $\omega_{k}$ to denote the $k^{th}$ element
of $\omega\in\Omega$. Then $(\Omega,d_{\Omega})$ is a compact metric space,
where the metric $d_{\Omega}$ is defined by $d_{\Omega}(\sigma,\omega)=0$ when
$\sigma=\omega$ and $d_{\Omega}(\sigma,\omega)=2^{-k}$ when $k$ is the least
index for which $\sigma_{k}\neq\omega_{k}$. We call $\Omega$ the code space
associated with the IFS\ $\mathcal{F}$.

Let $\sigma\in\Omega$ and $x\in\mathbb{X}$. Then, using the contractivity of
$\mathcal{F}$, it is straightfoward to prove that
\[
\phi_{\mathcal{F}}(\sigma):=\lim_{k\rightarrow\infty}f_{\sigma_{1}}\circ
f_{\sigma_{2}}\circ...f_{\sigma_{k}}(x)
\]
exists, uniformly for $x$ in any fixed compact subset of $\mathbb{X}$, and
depends continuously on $\sigma$. See for example \cite{BaDe}, Theorem 3. Let
\[
A_{\mathcal{F}}=\{\phi_{\mathcal{F}}(\sigma):\sigma\in\Omega\}\text{.}%
\]
Then $A_{\mathcal{F}}\subset\mathbb{X}$ is called the attractor of
$\mathcal{F}$. The continuous function
\[
\phi_{\mathcal{F}}:\Omega\rightarrow A_{\mathcal{F}}%
\]
is called the address function of $\mathcal{F}$. We call $\phi_{\mathcal{F}%
}^{-1}(\{x\}):=\{\sigma\in\Omega:\phi_{\mathcal{F}}(\sigma)=x\}$ the set of
addresses of the point $x\in A_{\mathcal{F}}$.

Clearly $A_{\mathcal{F}}$ is compact, nonempty, and has the property
\[
A_{\mathcal{F}}=f_{1}(A_{\mathcal{F}})\cup f_{2}(A_{\mathcal{F}})\cup...\cup
f_{N}(A_{\mathcal{F}})\text{.}%
\]
Indeed, if we define $\mathbb{H(X)}$ to be the set of nonempty compact subsets
of $\mathbb{X}$, and we define $\mathcal{F}:\mathbb{H(X)\rightarrow H(X)}$ by%
\begin{equation}
\mathcal{F}(S)=f_{1}(S)\cup f_{2}(S)\cup...\cup f_{N}(S)\text{,}
\label{fseteq}%
\end{equation}
for all $S\in\mathbb{H(X)}$, then $A_{\mathcal{F}}$ can be characterized as
the unique fixed point of $\mathcal{F}$, see \cite{hutchinson}, section 3.2,
and \cite{williams}.

\begin{figure}
[ptb]
\begin{center}
\includegraphics[
height=2.0237in,
width=4.8845in
]%
{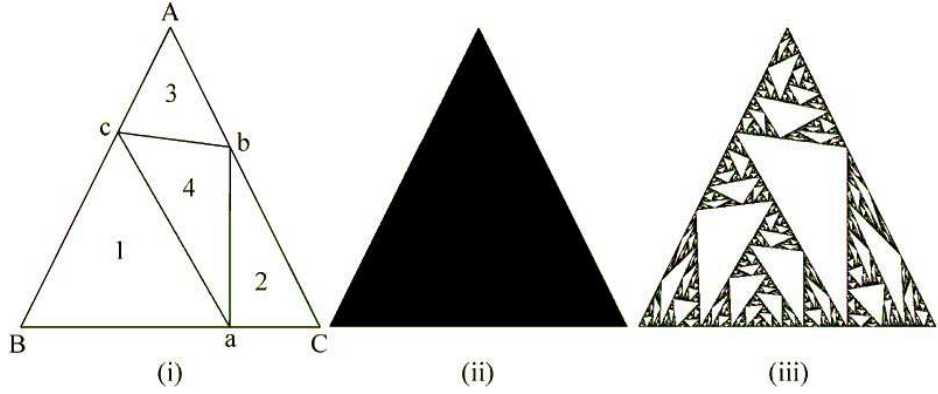}%
\caption{(i) The points used to define the affine transformations
$h_{n}:\mathbb{R}^{2}\rightarrow\mathbb{R}^{2}$ for $n=1,2,3,4$; (ii) sketch
of the attractor of the IFS $\{\mathbb{R}^{2};h_{1},h_{2},h_{3},h_{4}\}$;
(iii) sketch of the attractor of the IFS $\{\mathbb{R}^{2};h_{1},h_{2}%
,h_{3}\}$. Here $\alpha=0.65$, $\beta=0.3$, and $\gamma=0.4$.}%
\label{homtrisel}%
\end{center}
\end{figure}

IFSs may be used to represent diverse subsets of $\mathbb{R}^{2}$. For
example,\ let $A$, $B$, and $C$, denote three noncollinear points in
$\mathbb{R}^{2}$. Let $a$ denote a point on the line segment $AB$, let $b$
denote a point on the line segment $BC$ and let $c$ denote a point on the line
segment $CA$, such that $\{a,b,c\}\cap\{A,B,C\}=\varnothing$, see panel (i) of
Figure \ref{homtrisel}.%
Let $h_{1}:\mathbb{R}^{2}\rightarrow\mathbb{R}^{2}$ denote the unique affine
transformation such that
\[
h_{1}(ABC)=aBc\text{,}%
\]
by which we mean that $h_{1}$ maps $A$ to $a$, $B$ to $B$, and $C$ to $c$.
Using the same notation, let affine transformations $h_{2}$, $h_{3}$, and
$h_{4}$ be uniquely defined by%
\[
h_{2}(ABC)=abC\text{, }h_{3}(ABC)=Abc\text{, and }h_{4}(ABC)=abc\text{.}%
\]
Let $\mathcal{F}_{\alpha,\beta,\gamma}=\{\mathbb{R}^{2};h_{1},h_{2}%
,h_{3},h_{4}\}$ where $\alpha=|Bc|/|AB|,\beta=|Ca|/|BC|$ and $\gamma
=|Ab|/|CA|$. The attractor of $\mathcal{F}_{\alpha,\beta,\gamma}$ is the
filled triangle with vertices at $A$, $B$, and $C$. The attractor of the IFS
$\{\mathbb{R}^{2};h_{1},h_{2},h_{3}\}$ is an affine Sierpinski triangle, as
illustrated in (iii) in Figure \ref{homtrisel}.

For reference we note that when $A=(0,0)$, $B=(0,1),$ and $C=(0.5,1),$ the
transformations of the IFS $\mathcal{F}_{\alpha,\beta,\gamma}$ are given by
\[
h_{n}(x,y)=(a_{n}x+b_{n}y+c_{n},d_{n}x+e_{n}y+l_{n})
\]
with the parameters specified in Table \ref{ProjTable}. We will write
$\blacktriangle$ to denote the filled triangle $ABC$.%

\begin{table}[tbp] \centering
\begin{tabular}
[c]{|c|c|c|c|c|c|c|}\hline
$n$ & $a_{n}$ & $b_{n}$ & $c_{n}$ & $d_{n}$ & $e_{n}$ & $l_{n}$\\\hline
$1$ & $-1+\beta$ & $-\frac{1}{2}+\frac{1}{2}\beta+\frac{1}{2}\alpha$ &
$1-\beta$ & $0$ & $\alpha$ & $0$\\
$2$ & $\beta+\frac{1}{2}\gamma-\frac{1}{2}$ & $\frac{1}{2}\beta-\frac{1}%
{4}\gamma+\frac{1}{4}$ & $1-\beta$ & $1-\gamma$ & $\frac{1}{2}\gamma-\frac
{1}{2}$ & $0$\\
$3$ & $\frac{1}{2}\gamma$ & $-\frac{1}{2}+\frac{1}{2}\alpha-\frac{1}{4}\gamma$
& $\frac{1}{2}$ & $-\gamma$ & $-1+\alpha+\frac{1}{2}\gamma$ & $1$\\
$4$ & $\beta+\frac{1}{2}\gamma-\frac{1}{2}$ & $-\frac{3}{4}+\frac{1}{2}%
\beta+\frac{1}{2}\alpha-\frac{1}{4}\gamma$ & $1-\beta$ & $1-\gamma$ &
$\alpha-\frac{1}{2}+\frac{1}{2}\gamma$ & $0$\\\hline
\end{tabular}
\caption{ \label{ProjTable}}%
\end{table}%

\section{Chaos game}

When the underlying space is the euclidean plane, one way to sketch the
attactor $A_{\mathcal{F}}$ of an IFS $\mathcal{F}$ is to plot the set of
points%
\[
\widetilde{A}_{\mathcal{F}}=\{f_{\sigma_{1}}\circ f_{\sigma_{2}}%
\circ...f_{\sigma_{K}}(x):\sigma_{k}\in\{1,2,...,N\},k=1,2,...K\},
\]
for some $x\in\mathbb{X}$ and some integer $K$. The Hausdorff distance between
$A_{\mathcal{F}}$ and $\widetilde{A}_{\mathcal{F}}$ is bounded above by
$C\cdot l^{K}$ where the constant $C$ depends only on $\mathcal{F}$ and $x$.

A more efficient method is by means of a type of Markov Chain Monte Carlo
algorithm which we refer to as the chaos game. Starting from any point
$(x_{0},y_{0})\in\mathbb{R}^{2}$, a sequence of a million or more points
$\{(x_{k},y_{k})\}_{k=0}^{K}$ is computed recursively; at the $k^{th}$
iteration one of the functions of $\mathcal{F}$ is chosen at random,
independently of all other choices, and applied to $(x_{k-1},y_{k-1})$ to
produce $(x_{k},y_{k})$ which is plotted when $k\geq100$. The result will be
usually a sketch of the attractor of the IFS, accurate to within viewing resolution.

The reason that the chaos game yields, almost always, a "picture" of the
attractor of an IFS depends on Birkhoff's ergodic theorem, see for example
\cite{vrscay1}. The scholarly history of the chaos game is discussed in
\cite{kaijser} and \cite{stenflo3}, and appears to begin in 1935 with the work
of Onicescu and Mihok, \cite{onicescu}. Mandelbrot used a version of it to
help compute pictures of certain Julia sets, \cite{mandelbrot} pp.196-199; it
was introduced to IFS theory and developed by the author and coworkers, see
for example \cite{BaDe},\ \cite{barnsley}, \cite{marcaBerg}, and
\cite{JEergodic}, where the relevant theorems and much discussion can be
found. Its applications to fractal geometry were popularized initially by the
author and others, see for example \cite{BaDe}, \cite{devaney}, \cite{peak},
and \cite{peitgen}.

The sketches in panels (ii) and (iii) of Figure \ref{homtrisel} were computed
using the chaos game. At each iteration the function $h_{n}$ was selected with
probability proportional to the area of the triangle $h_{n}(ABC)$, for
$n=1,2,3,4$.

In section \ref{topsec} we show how the chaos game may be modified to
calculate examples of the fractal transformations that are the subject of this
article. Hopefully you will be inspired to try this new application of the
chaos game.

\section{\label{topsfnsec}The tops function}

We order the elements of $\Omega$ according to
\[
\sigma<\omega\text{ iff }\sigma_{k}>\omega_{k}%
\]
where $k$ is the least index for which $\sigma_{k}\neq\omega_{k}$. This is a
linear ordering, sometimes called the lexicographic ordering.

Notice that all elements of $\Omega$ are less than or equal to $\overline
{1}=11111...$ and greater than or equal to $\overline{N}=NNNNN....$. Also, any
pair of distinct elements of $\Omega$ is such that one member of the pair is
strictly greater than the other. In particular, the set of addresses of a
point $x\in A_{\mathcal{F}}$ is both closed and bounded above by $\overline
{1}$. It follows that $\phi_{\mathcal{F}}^{-1}(\{x\})$ possesses a unique
largest element. We denote this element by $\tau_{\mathcal{F}}(x)$.

\begin{definition}
Let $\mathcal{F}$ be a hyperbolic IFS with attractor $A_{\mathcal{F}}$ and
address function $\phi_{\mathcal{F}}:\Omega\rightarrow A_{\mathcal{F}}$. Let
\[
\tau_{\mathcal{F}}(x)=\max\{\sigma\in\Omega:\phi_{\mathcal{F}}(\sigma
)=x\}\text{ for all }x\in A_{\mathcal{F}}\text{.}%
\]
Then
\[
\Omega_{\mathcal{F}}:=\{\tau_{\mathcal{F}}(x):x\in A_{\mathcal{F}}\}
\]
is called the tops code space and
\[
\tau_{\mathcal{F}}:A_{\mathcal{F}}\rightarrow\Omega_{\mathcal{F}}%
\]
is called the tops function, for the IFS $\mathcal{F}$.
\end{definition}

Notice that the tops function $\tau_{\mathcal{F}}$ is one-to-one. It provides
a right-hand inverse to the address function, according to%
\[
\phi_{\mathcal{F}}\circ\tau_{\mathcal{F}}=i_{A_{\mathcal{F}}}%
\]
where $i_{A_{\mathcal{F}}}$ denotes the identity function on $A_{\mathcal{F}}$
and $\circ$ denotes composition of functions. Let
\[
\Phi_{\mathcal{F}}:\Omega_{\mathcal{F}}\rightarrow A_{\mathcal{F}}%
\]
denote the restriction of $\phi_{\mathcal{F}}$ to $\Omega_{\mathcal{F}}$,
defined by $\Phi_{\mathcal{F}}(\sigma)=\phi_{\mathcal{F}}(\sigma)$ for all
$\sigma\in\Omega_{\mathcal{F}}$. Then $\Phi_{\mathcal{F}}$ is the inverse of
$\tau_{\mathcal{F}}$, namely%
\[
\Phi_{\mathcal{F}}=\tau_{\mathcal{F}}^{-1}\text{.}%
\]
We note that although $\Phi_{\mathcal{F}}$ is one-to-one, onto, and
continuous, $\tau_{\mathcal{F}}$ may not be continuous. Let $\overline{\Omega
}_{\mathcal{F}}$ denote the closure of $\Omega_{\mathcal{F}}$, treated as a
subset of the metric space $(\Omega,d_{\Omega})$. Let
\[
\overline{\Phi}_{\mathcal{F}}:\overline{\Omega}_{\mathcal{F}}\rightarrow
A_{\mathcal{F}}%
\]
denote the restriction of $\phi_{\mathcal{F}}$ to $\overline{\Omega
}_{\mathcal{F}}$. Then $\overline{\Phi}_{\mathcal{F}}$ is continuous and onto.
Notice that the ranges of $\overline{\Phi}_{\mathcal{F}}$ and $\Phi
_{\mathcal{F}}$ are both equal to $A_{\mathcal{F}}$ because $A_{\mathcal{F}}$
is closed.

\section{Fractal transformations}

Let $\mathcal{G}$ denote a hyperbolic IFS that consists of $N$ functions. Then
$\phi_{\mathcal{G}}\circ\tau_{\mathcal{F}}:A_{\mathcal{F}}\rightarrow
A_{\mathcal{G}}$ is a mapping from the attractor of $\mathcal{F}$ into the
attractor of $\mathcal{G}$. We refer to $\phi_{\mathcal{G}}\circ
\tau_{\mathcal{F}}$ as a fractal transformation.

In order to illustrate transformations between subsets of $\mathbb{R}^{2}$ we
use pictures. We define a picture function to be a function of the form
\[
\mathfrak{P}:D_{\mathfrak{P}}\subset\mathbb{R}^{2}\rightarrow\mathfrak{C}%
\]
where $\mathfrak{C}$ is a color space. A picture function $\mathfrak{P}$
assigns a unique color to each point in its domain $D_{\mathfrak{P}}$. For
example we may have $\mathfrak{C=}\{0,1,...255\}^{3}$ and each point of
$\mathfrak{C}$ may specify the red, green, and blue components of a color. A
picture in the non-mathematical sense may be thought of as a physical
representation of the graph of a picture function.

If $T:D_{\mathfrak{Q}}\subset\mathbb{R}^{2}\rightarrow D_{\mathfrak{P}}$ then
$\mathfrak{Q}=\mathfrak{P}\circ T$ denotes a picture whose domain is
$D_{\mathfrak{Q}}$. We can obtain insights into the nature of $T$ by comparing
the picture functions $\mathfrak{P}$ and $\mathfrak{P}\circ T$, where
$\mathfrak{P}$ represents a given picture which may be varied. We will use
this method to illustrate fractal transformations.

\begin{figure}
[ptb]
\begin{center}
\includegraphics[
height=3.2085in,
width=4.8075in
]%
{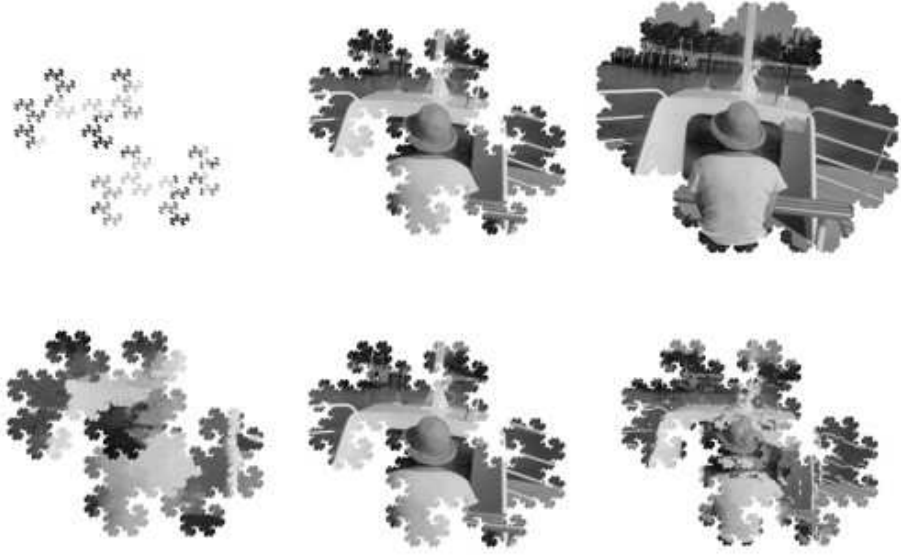}%
\caption{Examples of fractal transformations are illustrated using the three
picture functions $\mathfrak{P}_{\mathcal{G}}$, $\mathfrak{P}_{\mathcal{F}}$,
and $\mathfrak{P}_{\mathcal{H}}$ shown in the top row, from left to right. The
domains of these functions are the attractors $A_{\mathcal{G}}$,
$A_{\mathcal{F}}$, $A_{\mathcal{H}}$ of three IFSs $\mathcal{F}$,
$\mathcal{G}$, and $\mathcal{H}$, defined in the text. The bottom row
illustrates the pictures $\mathfrak{P}_{\mathcal{G}}\mathfrak{\circ}%
\phi_{\mathcal{G}}\circ\tau_{\mathcal{F}},$ $\mathfrak{P}_{\mathcal{F}%
}\mathfrak{\circ}\phi_{\mathcal{F}}\circ\tau_{\mathcal{F}},$ and
$\mathfrak{P}_{\mathcal{H}}\mathfrak{\circ}\phi_{\mathcal{H}}\circ
\tau_{\mathcal{F}}$.}%
\label{boatdragonsgray}%
\end{center}
\end{figure}

Let
\[
\mathcal{F}=\{\mathbb{C};f_{1}(z)=sz-1,f_{2}(z)=sz+1,\text{ for all }%
z\in\mathbb{C}\},
\]%
\[
\mathcal{G}=\{\mathbb{C};g(z)=s_{1}z-1,g(z)=s_{1}z+1,\text{ for all }%
z\in\mathbb{C}\},
\]%
\[
\mathcal{H}=\{\mathbb{C};h(z)=s_{2}z-1,h(z)=s_{2}z+1,\text{ for all }%
z\in\mathbb{C}\},
\]
where $\mathbb{C}$, denotes the complex plane, $s=0.5(1+i)$, $s_{1}%
=0.44(1+i)$, and $s_{2}=0.535(1+i)$. We denote the attractors of these IFSs by
$A_{\mathcal{F}}$, $A_{\mathcal{G}}$, and $A_{\mathcal{H}}$. In the top row of
Figure \ref{boatdragonsgray} we illustrate, from left to right, the three
picture functions, $\mathfrak{P}_{\mathcal{G}}:A_{\mathcal{G}}\rightarrow
\mathfrak{C}$, $\mathfrak{P}_{\mathcal{F}}:A_{\mathcal{F}}\rightarrow
\mathfrak{C}$, and $\mathfrak{P}_{\mathcal{H}}:A_{\mathcal{H}}\rightarrow
\mathfrak{C}$. These pictures were obtained by masking a single original
digital picture, whose domain we took to be $\{z=x+iy\in\mathbb{C}:$ $-3.5\leq
x\leq3.5,$ $-3.5\leq y\leq3.5\}$, by the complement of each of the sets
$A_{\mathcal{G}}$, $A_{\mathcal{F}}$, and $A_{\mathcal{H}}$.

The attractor $A_{\mathcal{F}}$ is a so-called twin-dragon fractal. It is an
example of a just-touching attractor: that is, $f_{1}(A_{\mathcal{F}})\cap
f_{2}(A_{\mathcal{F}})$ is non-empty and equals $f_{1}(\partial A_{\mathcal{F}%
})\cap f_{2}(\partial A_{\mathcal{F}})$ where $\partial A_{\mathcal{F}}$
denotes the boundary of $A_{\mathcal{F}}$. This contrasts with $A_{\mathcal{G}%
}$ which is totally disconnected, perfect, and in fact homeomorphic to the
classical cantor set. This also contrasts with $A_{\mathcal{H}}$ which is such
that there exists a disk in $\mathbb{R}^{2}$, of non-zero radius, which is
contained in $h_{1}(A_{\mathcal{H}})\cap h_{2}(A_{\mathcal{H}})$.

The bottom row of Figure \ref{boatdragonsgray} illustrates the pictures, from
left to right, $\mathfrak{P}_{\mathcal{G}}\mathfrak{\circ}\phi_{\mathcal{G}%
}\circ\tau_{\mathcal{F}},$ $\mathfrak{P}_{\mathcal{F}}\mathfrak{\circ}%
\phi_{\mathcal{F}}\circ\tau_{\mathcal{F}},$ and $\mathfrak{P}_{\mathcal{H}%
}\mathfrak{\circ}\phi_{\mathcal{H}}\circ\tau_{\mathcal{F}}$. They were
computed by a variant of the chaos game as explained in section \ref{topsec}.
The domain of each of these pictures is $A_{\mathcal{F}}$. We notice that
$\mathfrak{P}_{\mathcal{F}}\mathfrak{\circ}\phi_{\mathcal{F}}\circ
\tau_{\mathcal{F}}=\mathfrak{P}_{\mathcal{F}}$, which is true regardless of
the choice of IFS $\mathcal{F}$ since $\phi_{\mathcal{F}}\circ\tau
_{\mathcal{F}}$ is the identity on $A_{\mathcal{F}}$. We notice that both
$\mathfrak{P}_{\mathcal{G}}\mathfrak{\circ}\phi_{\mathcal{G}}\circ
\tau_{\mathcal{F}}$ and $\mathfrak{P}_{\mathcal{H}}\mathfrak{\circ}%
\phi_{\mathcal{H}}\circ\tau_{\mathcal{F}}$ have features in common with the
underlying digital picture; for example, $\mathfrak{P}_{\mathcal{G}%
}\mathfrak{\circ}\phi_{\mathcal{G}}\circ\tau_{\mathcal{F}}$ displays something
like the texture of the hat, near the middle of the bottom-left image. The
bottom right image shows parts of the hat, repeated several times, and some
clearly delineated small twin-dragon tiles.%

We are led to consider the following questions. Under what conditions on
general IFSs $\mathcal{F}$ and $\mathcal{G}$ is the fractal transformation
$\phi_{\mathcal{G}}\circ\tau_{\mathcal{F}}$ continuous? When does it provide a
homeomorphism between $A_{\mathcal{F}}$ and $A_{\mathcal{G}}$?

\begin{definition}
The address structure of\textbf{ }$\mathcal{F}$ is defined to be the set of
sets
\[
\mathcal{C}_{\mathcal{F}}=\{\phi_{\mathcal{F}}^{-1}(\{x\})\cap\overline
{\Omega}_{\mathcal{F}}:x\in A_{\mathcal{F}}\}\text{.}%
\]

\end{definition}

The address structure of an IFS is a certain partition of $\overline{\Omega
}_{\mathcal{F}}$. Let $\mathcal{C}_{\mathcal{G}}$ denote the address structure
of $\mathcal{G}$. Let us write $\mathcal{C}_{\mathcal{F}}$ $\prec
\mathcal{C}_{\mathcal{G}}$ to mean that for each\ $S\in\mathcal{C}%
_{\mathcal{F}}$ there is $T\in\mathcal{C}_{\mathcal{G}}$ such that $S\subset
T$. Notice that if $\mathcal{C}_{\mathcal{F}}=\mathcal{C}_{\mathcal{G}}$ then
$\Omega_{\mathcal{F}}=$ $\Omega_{\mathcal{G}}$. Some examples of address
structures are given in section \ref{AddressSec}.

\begin{theorem}
Let $\mathcal{F}$ and $\mathcal{G}$ be two hyperbolic IFSs such that
$\mathcal{C}_{\mathcal{F}}\prec\mathcal{C}_{\mathcal{G}}$. Then the fractal
transformation $\phi_{\mathcal{G}}\circ\tau_{\mathcal{F}}:A_{\mathcal{F}%
}\rightarrow A_{\mathcal{G}}$ is continuous. If $\mathcal{C}_{\mathcal{F}%
}=\mathcal{C}_{\mathcal{G}}$ then $\phi_{\mathcal{G}}\circ\tau_{\mathcal{F}}$
is a homeomorphism.
\end{theorem}

The proof relies on a standard result in topology, Lemma 2 below, which we
present in the context of metric spaces.

\begin{lemma}
(cf. \cite{Mendelson}, bottom of p.194.) Let $F:X\rightarrow Y$ be a
continuous mapping from a compact metric space $X$ onto a metric space $Y$.
Then $S\subset Y$ is open if and only if $F^{-1}(S)\subset X$ is open.
\end{lemma}

\begin{proof}
If $S\subset Y$ is open then $F^{-1}(S)\subset X$ is open because
$F:X\rightarrow Y$ is continuous. Suppose that $F^{-1}(S)$ is open. Then
$X\backslash F^{-1}(S)$ is closed. But a closed subset of a compact metric
space is compact. The continuity of $F$ now implies that $F(X\backslash
F^{-1}(S))$ is compact and hence closed. But $F(X\backslash F^{-1}%
(S))=Y\backslash S$. Hence $S$ is open.
\end{proof}

\begin{lemma}
(cf. \cite{Mendelson}, Proposition 7.4 on p.195.) Let $F:X\rightarrow Y$ be a
continuous mapping from a compact metric space $X$ onto a metric space $Y$.
Let $H:Y\rightarrow Z$ where $Z$ is a metric space.\ Let $H\circ
F:X\rightarrow Z$ be continuous. Then $H:Y\rightarrow Z$ is continuous.
\end{lemma}

\begin{proof}
Let $O\subset Z$ be open. Then $(H\circ F)^{-1}(O)=F^{-1}(H^{-1}(O))$ is open.
But then by Lemma 1 $H^{-1}(O)$ is open. Hence $H:Y\rightarrow Z$ is continuous.
\end{proof}

\begin{proof}
[Proof of Theorem 2]In Lemma 2 we set $X=\overline{\Omega}_{\mathcal{F}}$,
$Y=A_{\mathcal{F}}$, and $Z=A_{\mathcal{G}}$. We choose $F:X\rightarrow Y$ to
be $\overline{\Phi}_{\mathcal{F}}:\overline{\Omega}_{\mathcal{F}}\rightarrow
A_{\mathcal{F}}$. Then $F:X\rightarrow Y$ is a continuous mapping from a
compact metric space $X$ onto a metric space $Y$. We also choose
$H:Y\rightarrow Z$ to be
\[
H=\phi_{\mathcal{G}}\circ\tau_{\mathcal{F}}:A_{\mathcal{F}}\rightarrow
A_{\mathcal{G}}\text{.}%
\]
Now look at the function
\[
G:=H\circ F=\phi_{\mathcal{G}}\circ\tau_{\mathcal{F}}\circ\overline{\Phi
}_{\mathcal{F}}:\overline{\Omega}_{\mathcal{F}}\rightarrow A_{\mathcal{G}%
}\text{.}%
\]
If $\sigma\in\overline{\Omega}_{\mathcal{F}}$, then both $\sigma$ and
$(\tau_{\mathcal{F}}\circ\overline{\Phi}_{\mathcal{F}})(\sigma)$ belong to the
same set in the $\mathcal{C}_{\mathcal{F}}$. Since $\mathcal{C}_{\mathcal{F}%
}\prec\mathcal{C}_{\mathcal{G}}$ it follows that both $\sigma$ and
$(\tau_{\mathcal{F}}\circ\overline{\Phi}_{\mathcal{F}})(\sigma)$ belong to the
same set in $\mathcal{C}_{\mathcal{G}}$. It follows that%
\[
(\phi_{\mathcal{G}}\circ\tau_{\mathcal{F}}\circ\overline{\Phi}_{\mathcal{F}%
})(\sigma)=\phi_{\mathcal{G}}(\sigma)\text{ for all }\sigma\in\overline
{\Omega}_{\mathcal{F}}\text{.}%
\]
But $\phi_{\mathcal{G}}:\Omega\rightarrow A_{\mathcal{G}}$ is continuous.
Hence $G$ is continuous.

We have shown that the conditions in Lemma 2 hold. It follows that
$H=\phi_{\mathcal{G}}\circ\tau_{\mathcal{F}}$ is continuous.

When $\mathcal{C}_{\mathcal{F}}=\mathcal{C}_{\mathcal{G}}$ it is readily
verified that $\phi_{\mathcal{G}}\circ\tau_{\mathcal{F}}:A_{\mathcal{F}%
}\rightarrow A_{\mathcal{G}}$ is one-to-one and onto and that its inverse is
$\phi_{\mathcal{F}}\circ\tau_{\mathcal{G}}$. Also $\mathcal{C}_{\mathcal{F}%
}=\mathcal{C}_{\mathcal{G}}$ implies $\mathcal{C}_{\mathcal{G}}\prec
\mathcal{C}_{\mathcal{F}}$ and so, by the first part of the theorem,
$\phi_{\mathcal{F}}\circ\tau_{\mathcal{G}}$ is continuous. Hence
$\phi_{\mathcal{G}}\circ\tau_{\mathcal{F}}:A_{\mathcal{F}}\rightarrow
A_{\mathcal{G}}$ is a homeomorphism.
\end{proof}

%

\begin{figure}
[ptb]
\begin{center}
\includegraphics[
height=2.1465in,
width=4.8075in
]%
{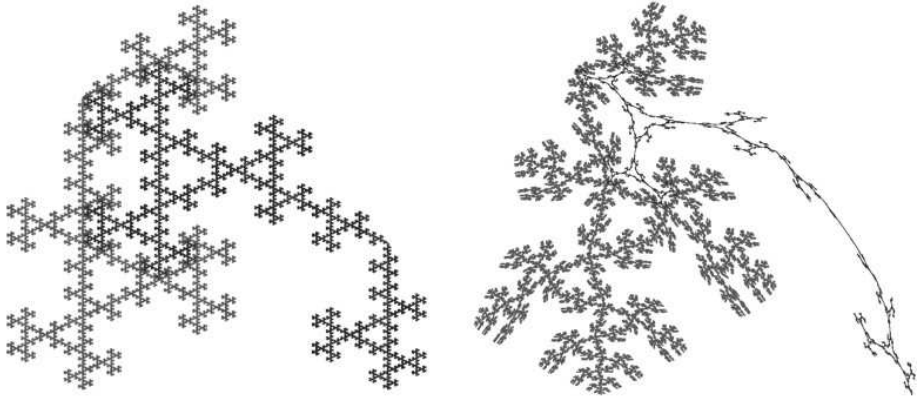}%
\caption{Before, on the left, and after, on the right, a fractal
homeomorphism. See text.}%
\label{trihombefoandaft2}%
\end{center}
\end{figure}

\section{\label{AddressSec}Examples of address structures}

\subsection{Backwards orbits}

Let $x\in\mathcal{A}_{\mathcal{F}}$. Let $\sigma\in\Omega$ be such that
$\phi_{\mathcal{F}}(\sigma)=x$. Assume that the $f_{n}$'s are one-to-one. Then
define $x_{0}=x$ and $x_{k}=f_{\sigma_{k}}^{-1}(x_{k-1})$ for $k=1,2,3,...$.
Notice that $x_{k}=\phi_{\mathcal{F}}(\sigma_{k}\sigma_{k+1}\sigma_{k+2}...)$.
We call $\{x_{k}\}_{k=0}^{\infty}$ a backwards orbit of $x$ (under the IFS
$\mathcal{F)}$.

The set of all addresses of $x$ can be calculated by following all possible
backwards orbits of $x$. Define a sequence of points $\{\widetilde{x}%
_{k}\}_{k=0}^{\infty}$ in $A_{\mathcal{F}}$ and an address $\widetilde{\sigma
}=\widetilde{\sigma}_{1}\widetilde{\sigma}_{2}\widetilde{\sigma}_{3}%
...\in\Omega$, as follows. Let $\widetilde{x}_{0}=x$. For each $k=1,2,3,...$
first choose%
\[
\widetilde{\sigma}_{k}\in\{n\in\{1,2,...,N\}:\widetilde{x}_{k-1}\in
f_{n}(A_{\mathcal{F}})\}
\]
and then define%
\[
\widetilde{x}_{k}=f_{\widetilde{\sigma}_{k}}^{-1}(\widetilde{x}_{k-1})\text{.}%
\]
Then $\widetilde{\sigma}\in\phi_{\mathcal{F}}^{-1}(\{x\})$ and all
$\widetilde{\sigma}\in\phi_{\mathcal{F}}^{-1}(\{x\})$ can be obtained in this manner.

\subsection{Some notation}

We use the notation $(PQ)=PQ\backslash\{P,Q\}$ to denote the straight line
segment which connects the two points $P$ and $Q$ in $\mathbb{R}^{2},$ without
its endpoints. We write $\Omega^{\prime}$ to denote the set of all finite
length strings of symbols from the alphabet $\{1,2,...,N\}$, including the
empty string "$\varnothing$" . We write $\left\vert \sigma\right\vert $ to
denote the length of $\sigma\in\Omega^{\prime}$. We define $\omega
\sigma=\omega_{1}\omega_{2}...\omega_{\left\vert \omega\right\vert }\sigma
_{1}\sigma_{2}...\sigma_{\left\vert \sigma\right\vert }$ for all
$\omega,\sigma\in\Omega^{\prime}$. Similarly we define $\omega\sigma
=\omega_{1}\omega_{2}...\omega_{\left\vert \omega\right\vert }\sigma_{1}%
\sigma_{2}...$ for all $\omega\in\Omega^{\prime},\sigma\in\Omega$. We write
$S:\Omega^{\prime}\rightarrow\Omega^{\prime}$ to denote the shift operator
defined by $S(\sigma)=\sigma_{2}\sigma_{3}...\sigma_{\left\vert \sigma
\right\vert }$ when $\left\vert \sigma\right\vert \geq1$ and $S($%
"$\varnothing$"$)=$"$\varnothing$". We write $f_{\sigma}=f_{\sigma_{1}}\circ
f_{\sigma_{2}}\circ...\circ f_{\sigma_{\left\vert \sigma\right\vert }}$ for
all $\sigma\in\Omega^{\prime}$ with $\left\vert \sigma\right\vert \geq1$ and
$f_{"\varnothing"}$ denotes the identity function.

\subsection{\label{example1}Example 1}

An interesting example of address structures is provided by the IFS
$\mathcal{F}=\mathcal{F}_{\alpha,\beta,\gamma}=\{\mathbb{R}^{2};h_{1}%
,h_{2},h_{3},h_{4}\}$, introduced at the end of section \ref{hypersec}. Here
we prove that
\begin{equation}
\mathcal{C}_{\mathcal{F}_{\alpha,\beta,\gamma}}=\mathcal{C}_{\mathcal{F}%
_{\widetilde{\alpha},\widetilde{\beta},\widetilde{\gamma}}}
\label{equalstructs}%
\end{equation}
for all $\alpha,\beta,\gamma,\widetilde{\alpha},\widetilde{\beta}%
,\widetilde{\gamma}\in(0,1)$ by calculating the address structure
$\mathcal{C}_{\mathcal{F}_{\alpha,\beta,\gamma}}$.

In this case there is only one backwards orbit $\{x_{k}\}_{k=0}^{\infty}$ of
each $x\in\mathcal{A}_{\mathcal{F}}$. This is because of the form of
$\mathcal{F}$, and because the $h_{n}$'s are affine and so preserve ratios of
distances between points which lie on any given straight line: for example if
$x\in ab=h_{2}(\blacktriangle)\cap h_{4}(\blacktriangle)$ then $h_{2}%
^{-1}(x)=h_{4}^{-1}(x)$. Indeed, the mapping $T:\blacktriangle\rightarrow
\blacktriangle$ defined as in equations \ref{partition} is continuous and the
backwards orbit of $x$ is the same as the orbit of $x$ under $T$ treated as a
dynamical system.

Any point $x_{K}\in\{x_{k}\}_{k=0}^{\infty}$ on the backwards orbit of $x,$
such that more than one map $h_{n}^{-1}$ may be applied, is such that $x_{K}$
belongs to the set
\[%
{\textstyle\bigcup\limits_{i\neq j}}
h_{i}(\blacktriangle)\cap h_{j}(\blacktriangle)=(ab)\cup(bc)\cup
(ca)\cup\{a,b,c\}\text{.}%
\]
If $x_{K}\in(ab)$ then $\sigma_{K}\in\{2,4\}$, if $x_{K}=a$ then $\sigma
_{K}\in\{1,2,4\}$, and so on. For example, if $x=\phi_{\mathcal{F}}(\omega
_{1}\omega_{2}\omega_{3}...)$ and the only point on the backwards orbit of $x$
which lies in $(ab)\cup(bc)\cup(ca)\cup\{a,b,c\}$ is $x_{K}\in(ab)$ then
$\phi_{\mathcal{F}}^{-1}(\{x\})=\{\sigma\in\Omega:$ $\sigma_{k}=\omega_{k}$
for all $k\neq K$, and $\sigma_{K}\in\{2,4\}\}$.

Let $\triangle$ denote the boundary of $\blacktriangle$ as a subset of
$\mathbb{R}^{2}$, $\triangledown=(ab)\cup(bc)\cup(ca),$ and
\[
\Xi=\triangle\cup\triangledown\cup(%
{\textstyle\bigcup\limits_{\{\sigma\in\Omega^{\prime}:\left\vert
\sigma\right\vert \geq1\}}}
h_{\sigma}(\triangledown))\text{.}%
\]
Then

(i) each of the sets $\triangle,\triangledown,$ $h_{1}(\triangledown),$
$h_{2}(\triangledown),$ $h_{3}(\triangledown),$ $h_{4}(\triangledown),$
$h_{11}(\triangledown),$ $h_{12}(\triangledown),...$ is disjoint;

(ii) $T(\Xi)=\mathcal{F}(\Xi)=\Xi$;

(iii) $T$ is one-to-one on $h_{\sigma}(\triangledown)$ and $T(h_{\sigma
}(\triangledown))=h_{S(\sigma)}(\triangledown)$ for all $\sigma\in
\Omega^{\prime}$ with $\left\vert \sigma\right\vert \geq1$;

(iv) $T$ is one-to-one on $\triangledown$ and $T(\triangledown)=\triangle
\backslash\{A,B,C\}$;

(v) $T$ is two-to-one on $\triangle$ and $T(\triangle\backslash
\{A,B,C\})=T(\triangle)=\triangle$.

\begin{figure}
[ptb]
\begin{center}
\includegraphics[
height=4.6726in,
width=4.8845in
]%
{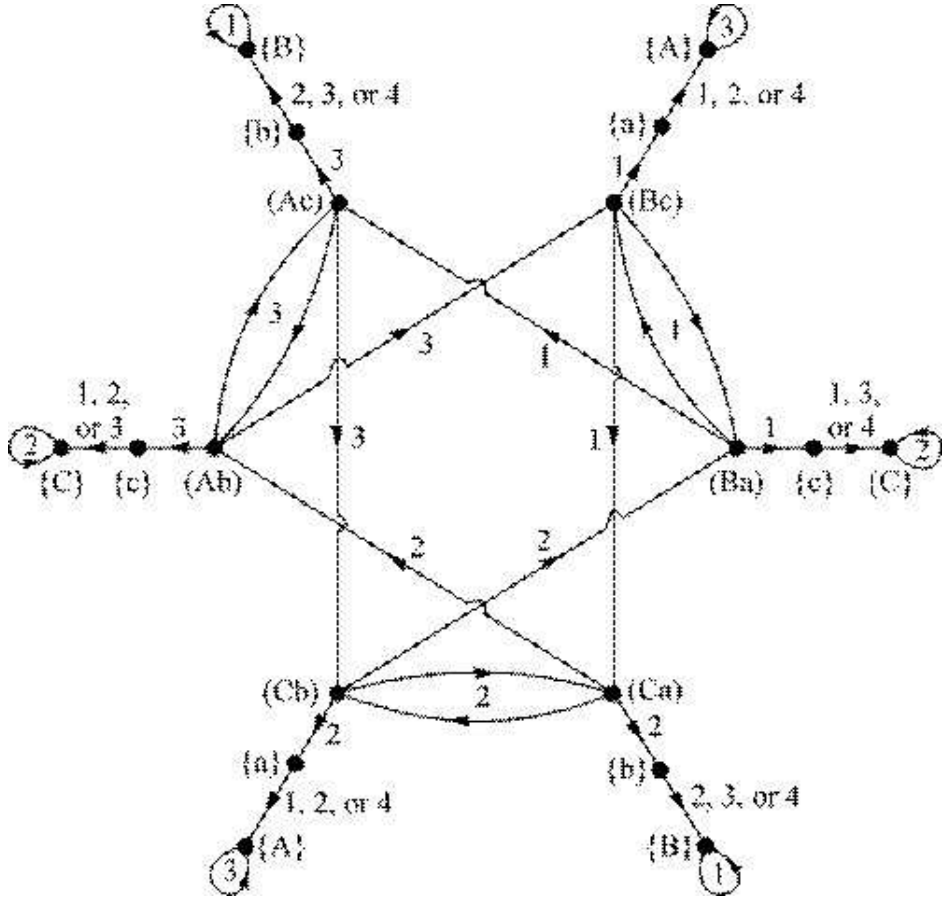}%
\caption{This graph can be used to compute the IFS addresses of all points
which lie on the boundary of the triangle $\blacktriangle$.}%
\label{markov}%
\end{center}
\end{figure}

The transformation $T$ maps $\triangle$ continuously onto itself. If $x$ goes
around $\triangle$ clockwise once, then $T(x)$ goes around $\triangle$
anticlockwise twice. It does so in such a way that $T(\{A\})=\{A\},$
$T((Ac))=(Ab)\cup\{b\}\cup(bC),$ $T(\{c\})=\{C\}$ and so on. This information
provides us with the directed graph, with labelled edges, shown in Figure
\ref{markov}. We denote this graph by $G$. It is such that that there is a
bijective correspondence between the points of $\triangle$ and the set of all
paths in $G$. A path in $G$ is obtained by starting at any node and
successively following edges in the directions of the arrows, yielding an
infinite sequences of edges. The set of addresses of the point represented by
a path in $G$ consists of all sequences of the numbers $\{1,2,3,4\}$ which can
be read off successivey from the path, with one symbol from each edge. For
example, the only possible address for $\{A\}$ is $\overline{3}=3333..$, and
the set of addresses of a point in $(Ab)$ may be $3\{1,2,$ or $3\}222...,$ or
$33\{2,3,$ or $4\}111..,.$ or one or more which begin $332222$.%

Now suppose $x\in h_{\sigma}(\triangledown)$ for some $\sigma\in\Omega
^{\prime}$ with $\left\vert \sigma\right\vert \geq1.$ Then by repeated
application of (iii) above we find that the first $\left\vert \sigma
\right\vert $ terms in any address of $x$ are precisely $\sigma$, and
$T^{\circ\left\vert \sigma\right\vert }(x)\in\triangledown$. Since
$T^{\circ\left\vert \sigma\right\vert }$ maps $h_{\sigma}(\triangledown)$
one-to-one onto $\triangledown$, it follows that the set of all sets of
addresses of all points in $h_{\sigma}(\triangledown)$ is the same as the set
of all sets of addresses of all points in $\triangledown$ after $\sigma$ has
been appended to the front of each of the latter addresses. So what is the set
of all sets of addresses of all points in $\triangledown$?

We have $\triangledown=(ab)\cup(bc)\cup(ca)$. Let us deal with $(ab)$. The
transformation $T$ maps $(ab)$ one-to-one onto $(AB)=(Ac)\cup\{c\}\cup(cB)$.
It follows that the set of all sets of addresses of all points in $(ab)$,
which we denote by $C(ab),$ is determined by the set of sets of addresses of
all points in $(Ac)\cup\{c\}\cup(cB)$, which we denote by $C((Ac)\cup
\{c\}\cup(cB))$. Specifically,
\[
C(ab)=\{\{\eta\sigma:\eta\in\{2,4\},\sigma\in\pi\}:\pi\in C((Ac)\cup
\{c\}\cup(cB))\}
\]
The set of sets of addresses $C((Ac)\cup\{c\}\cup(cB))$ corresponds to the set
of paths in $G$ which start at a node labelled $(Ac)$, $\{c\}$, or $(cB)$. We
can similarly describe the address structures of $(bc)$, and $(ca)$. We are
thus able, in principle, to write down the set of addresses of each $x\in$
$\Xi$; in particular, the set of all sets thus obtained does not depend on
$\alpha,\beta,$ or $\gamma$.

Next we deal with $\Xi^{C}=\blacktriangle\backslash\Xi$. Since all points with
multiple addresses lie in $\Xi$ and the backwards orbit of each point in
$\Xi^{C}$ lies in $\Xi^{C}$ it follows that each point in $\Xi^{C}$ has a
unique address, and in particular $\phi_{\mathcal{F}}^{-1}(\Xi^{C}%
)=\Omega\backslash\phi_{\mathcal{F}}^{-1}(\Xi)$ does not depend on
$\alpha,\beta,$ or $\gamma$.

Finally we note that $\overline{\phi_{\mathcal{F}}^{-1}(\Xi)}=\Omega$. Hence
$\Omega_{\mathcal{F}}=\Omega$. Hence $\Omega_{\mathcal{F}}=\phi_{\mathcal{F}%
}^{-1}(\Xi)\cup\phi_{\mathcal{F}}^{-1}(\Xi^{C}).$ Since the equivence class
structures of both $\phi_{\mathcal{F}}^{-1}(\Xi)$ and $\phi_{\mathcal{F}}%
^{-1}(\Xi^{C})$ do not depend on $\alpha,\beta,$ or $\gamma,$ it follows that
equation \ref{equalstructs} is true.

So, for example, let $\mathcal{F}=\mathcal{F}_{0.5,0.5,0.5}=\{\blacktriangle
;h_{1},h_{2},h_{3},h_{4}\}$ and $\mathcal{G=F}_{0.65,0.3,0.4}=\{\blacktriangle
;g_{1},g_{2},g_{3},g_{4}\}$. Then $\mathcal{C}_{\mathcal{F}}=\mathcal{C}%
_{\mathcal{G}}$ and, by Theorem 1, the fractal transformation $\phi
_{\mathcal{G}}\circ\tau_{\mathcal{F}}:\blacktriangle\rightarrow\blacktriangle$
is a homeomorphism. Figure \ref{trihombefoandaft2} illustrates the action of
this homeomorphism. The figure on the left shows the set $S$, defined to be
the union of the attractors of the two IFSs $\{\blacktriangle;h_{1}%
,h_{3},h_{4}\}$ and $\{\blacktriangle;h_{2},h_{3},h_{4}\}$. The image on the
right shows the set $\widetilde{S}$, defined to be the union of the attractors
of the IFSs $\{\blacktriangle;g_{1},g_{3},g_{4}\}$ and $\{\blacktriangle
;g_{2},g_{3},g_{4}\}$. The two sets are related by $\phi_{\mathcal{G}}%
\circ\tau_{\mathcal{F}}(S)=\widetilde{S}$. Figure \ref{threebirds} illustrates
two other homeomorphisms associated with the family $\mathcal{F}_{\alpha
,\beta,\gamma}$. These examples were computed as described in section
\ref{topsec}.

\subsection{\label{example2}Example 2}

\begin{figure}
[ptb]
\begin{center}
\includegraphics[
height=2.5356in,
width=4.8075in
]%
{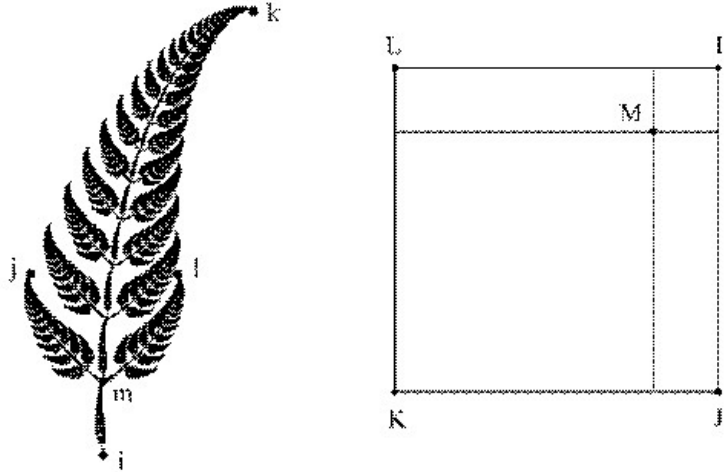}%
\caption{(i) Shows the points $i,j,k,l,m$ and (ii) shows the points
$I,J,K,L,M$. }%
\label{fernmaps}%
\end{center}
\end{figure}

An example of address structures $\mathcal{C}_{\mathcal{F}}$ and
$\mathcal{C}_{\mathcal{G}}$ such that $\mathcal{C}_{\mathcal{F}}%
\prec\mathcal{C}_{\mathcal{G}}$ and $\mathcal{C}_{\mathcal{F}}\neq
\mathcal{C}_{\mathcal{G}}$ is provided by taking $\mathcal{F}=\{\square
;f_{1},f_{2},f_{3},f_{4}\}$ and $\mathcal{G}=\{\square;g_{1},g_{2},g_{3}%
,g_{4}\}$ to be the IFSs of affine maps specified in Tables \ref{FernTable}
and \ref{CtsTable} respectively. Here $\square\subset\mathbb{R}^{2}$ denotes
the filled square with vertices at $I=(1,1)$, $J=(1,0)$, $K=(0,0)$, $J=(0,1)$.
The attractor $A_{\mathcal{F}}$ of $\mathcal{F}$ is represented by the fern
image in Figure \ref{fernmaps}. The attractor $A_{\mathcal{G}}$ of
$\mathcal{G}$ is $\square$.

The transformations of $\mathcal{F}$ are such\ that
\begin{align}
f_{1}(i)  &  =m,f_{1}(k)=k,f_{2}(i)=i,f_{2}(k)=m,\label{Faddress}\\
f_{3}(i)  &  =m,f_{3}(k)=l,f_{4}(i)=m,f_{4}(k)=j,\nonumber
\end{align}
where the points $i,j,k,l,m\in A_{\mathcal{F}}$ are approximately as labelled
in Figure \ref{fernmaps}. Furthermore $f_{p}(A_{\mathcal{F}})\cap
f_{q}(A_{\mathcal{F}})=m$ whenever $p,q\in\{1,2,3,4\}$ with $p\neq q$. It is
readily deduced that $k=\phi_{\mathcal{F}}(\overline{1}),i=\phi_{\mathcal{F}%
}(\overline{2}),m=\phi_{\mathcal{F}}(1\overline{2})=\phi_{\mathcal{F}%
}(2\overline{1})=\phi_{\mathcal{F}}(3\overline{2})=\phi_{\mathcal{F}%
}(4\overline{2})$, that
\[
\Omega_{\mathcal{F}}=\{\sigma\in\Omega:S^{\circ n}(\sigma)\notin
\{2\overline{1},3\overline{2},4\overline{2}\}\text{ for all }n\in
\{0,1,2,...\}\},
\]
that $\overline{\Omega}_{\mathcal{F}}=\Omega$ and that the address structure
of $\mathcal{F}$ is
\[
\mathcal{C}_{\mathcal{F}}=\mathcal{C}_{\mathcal{F}}^{(1)}\cup\mathcal{C}%
_{\mathcal{F}}^{(2)}%
\]
where
\begin{align*}
\mathcal{C}_{\mathcal{F}}^{(1)}  &  =\{\{\sigma\}:\sigma\in\Omega,S^{\circ
n}(\sigma)\notin\{1\overline{2},2\overline{1},3\overline{2},4\overline
{2}\}\text{ for all }n\in\{0,1,2,...\}\},\\
\mathcal{C}_{\mathcal{F}}^{(2)}  &  =\{\{\sigma^{\prime}1\overline{2}%
,\sigma^{\prime}2\overline{1},\sigma^{\prime}3\overline{2},\sigma^{\prime
}4\overline{2}\}:\sigma^{\prime}\in\Omega^{\prime}\}.
\end{align*}

To determine the address structure of $\mathcal{G}$, we note that $\square$ is
the union of four rectangular tiles $g_{n}(\square)$ which share portions of
their boundaries. The transformations of $\mathcal{G}$ are such that
\begin{align}
g_{1}(I)  &  =M,g_{1}(K)=K,g_{2}(I)=I,g_{2}(K)=M,\label{Gaddress}\\
g_{3}(I)  &  =M,g_{3}(K)=L,g_{4}(I)=M,g_{4}(K)=J,\nonumber
\end{align}
where the points $I,J,K,L,M\in A_{\mathcal{G}}$ are approximately as labelled
in Figure \ref{fernmaps}.

Note that equations \ref{Gaddress} are the same as equations \ref{Faddress}
upon substitution of $f_{1},f_{2},f_{3},f_{4},i,j,k,l,$ and $m$, by
$g_{1},g_{2},g_{3},g_{4},I,J,K,L,$ and $M$ respectively. It is readily deduced
that $K=\phi_{\mathcal{G}}(\overline{1}),I=\phi_{\mathcal{G}}(\overline
{2}),M=\phi_{\mathcal{G}}(1\overline{2})=\phi_{\mathcal{G}}(2\overline
{1})=\phi_{\mathcal{G}}(3\overline{2})=\phi_{\mathcal{G}}(4\overline{2})$, and
that $\overline{\Omega}_{\mathcal{G}}=\Omega$. As a consequence $\mathcal{C}%
_{\mathcal{F}}\prec\mathcal{C}_{\mathcal{G}}$: if $s\in\mathcal{C}%
_{\mathcal{F}}$ then either $s\in\mathcal{C}_{\mathcal{F}}^{(1)}$ or
$s\in\mathcal{C}_{\mathcal{F}}^{(2)}$; if $s\in\mathcal{C}_{\mathcal{F}}%
^{(1)}$ then $s$ is a singleton and, since $\mathcal{C}_{\mathcal{G}}$ is a
partition of $\Omega$, there must be $t\in\mathcal{C}_{\mathcal{G}}$ such that
$s\subset t$; if $s\in\mathcal{C}_{\mathcal{F}}^{(2)}$ then $s=\{\sigma
^{\prime}1\overline{2},\sigma^{\prime}2\overline{1},\sigma^{\prime}%
3\overline{2},\sigma^{\prime}4\overline{2}\}$ for some $\sigma^{\prime}%
\in\Omega^{\prime}$, and since $M=\phi_{\mathcal{G}}(1\overline{2}%
)=\phi_{\mathcal{G}}(2\overline{1})=\phi_{\mathcal{G}}(3\overline{2}%
)=\phi_{\mathcal{G}}(4\overline{2}),$ it follows that $\mathcal{C}%
_{\mathcal{G}}$ contains a set that contains $s$. Hence $\mathcal{C}%
_{\mathcal{F}}\prec\mathcal{C}_{\mathcal{G}}$ and, by Theorem 1, the fractal
transformation $\phi_{\mathcal{G}}\circ\tau_{\mathcal{F}}$ from the
fern-shaped set onto $\square$ is continuous. This transformation is
illustrated in Figure \ref{ferns}, as described at the start of section
\ref{topsec}. Note however that in this case $\mathcal{C}_{\mathcal{F}}%
\neq\mathcal{C}_{\mathcal{G}}$ because there is a set in $\mathcal{C}%
_{\mathcal{G}}$ which consist of a pair of distinct addresses, whereas all
sets in $\mathcal{C}_{\mathcal{F}}$ contain either one or four distinct addresses.

If, in this example, we change $\mathcal{G}$ to $\widetilde{\mathcal{G}}$
specified in Table \ref{DiscTable} then the attractor is still the filled
square, that is $A_{\widetilde{\mathcal{G}}}=\square$, but Equation
\ref{Gaddress} no longer holds and we can show that the fractal transformation
$\phi_{\widetilde{\mathcal{G}}}\circ\tau_{\mathcal{F}}$ from the fern-shaped
set onto $\square$ is not continuous. This lack of continuity is illustrated
in Figure \ref{ferns}, as described in section \ref{topsec}.
%

\begin{table}[tbp] \centering
\begin{tabular}
[c]{|c|c|c|c|c|c|c|}\hline
$n$ & $a_{n}$ & $b_{n}$ & $c_{n}$ & $d_{n}$ & $e_{n}$ & $l_{n}$\\\hline
$1$ & $0.85$ & $-0.05$ & $0.125$ & $0.05$ & $0.85$ & $-0.039$\\
$2$ & $0.06$ & $0.02$ & $0.45$ & $0.0$ & $0.165$ & $0.835$\\
$3$ & $0.17$ & $0.22$ & $0.195$ & $-0.22$ & $0.17$ & $0.776$\\
$4$ & $-0.17$ & $-0.22$ & $0.805$ & $-0.22$ & $0.17$ & $0.776$\\\hline
\end{tabular}
\caption{ \label{FernTable}}%
\end{table}%
%

\begin{table}[tbp] \centering
\begin{tabular}
[c]{|c|c|c|c|c|c|c|}\hline
$n$ & $a_{n}$ & $b_{n}$ & $c_{n}$ & $d_{n}$ & $e_{n}$ & $l_{n}$\\\hline
$1$ & $0.8$ & $0.0$ & $0.0$ & $0.0$ & $0.8$ & $0.0$\\
$2$ & $0.2$ & $0.0$ & $0.8$ & $0.0$ & $0.8$ & $0.2$\\
$3$ & $-0.2$ & $0.0$ & $1.0$ & $0.0$ & $0.8$ & $0.0$\\
$4$ & $0.8$ & $0.0$ & $0.0$ & $0.0$ & $-0.2$ & $1.0$\\\hline
\end{tabular}
\caption{ \label{CtsTable}}%
\end{table}%
%

\begin{table}[tbp] \centering
\begin{tabular}
[c]{|c|c|c|c|c|c|c|}\hline
$n$ & $a_{n}$ & $b_{n}$ & $c_{n}$ & $d_{n}$ & $e_{n}$ & $l_{n}$\\\hline
$1$ & $-0.8$ & $0.0$ & $0.8$ & $0.0$ & $-0.8$ & $0.8$\\
$2$ & $-0.2$ & $0.0$ & $1.0$ & $0.0$ & $-0.2$ & $1.0$\\
$3$ & $0.8$ & $0.0$ & $0.0$ & $0.0$ & $0.2$ & $0.8$\\
$4$ & $0.2$ & $0.0$ & $0.8$ & $0.0$ & $0.8$ & $0.0$\\\hline
\end{tabular}
\caption{ \label{DiscTable}}%
\end{table}%

\section{\label{topsec}Pictures of tops functions}

When the underlying space is $\mathbb{R}^{2}$ we can use the chaos game to
compute illustrations of fractal transformations.

Let two\ hyperbolic IFSs
\[
\mathcal{F}:=\{\square;f_{1},...,f_{N}\}\text{ and }\mathcal{G}:=\{\square
;g_{1},...,g_{N}\}
\]
and a picture function%
\[
\mathfrak{P}:\square\rightarrow\mathfrak{C}%
\]
be given, where%
\[
\square:=\{(x,y)\in\mathbb{R}^{2}:0\leq x,y\leq1\}\text{.}%
\]
Let $\mathfrak{P}_{\mathcal{G}}$ denote $\mathfrak{P}$ restricted to
$A_{\mathcal{G}}$, that is $\mathfrak{P}_{\mathcal{G}}=\mathfrak{P}|_{A_{G}}$.
Then we define a new picture%
\[
\mathfrak{P}_{\mathcal{F}}:A_{\mathcal{F}}\rightarrow\mathfrak{C}%
\]
by
\[
\mathfrak{P}_{\mathcal{F}}=\mathfrak{P}_{\mathcal{G}}\circ\phi_{\mathcal{G}%
}\circ\tau_{\mathcal{F}}\text{.}%
\]
We say that $\mathfrak{P}_{\mathcal{F}}$ is defined by\textbf{\ }tops plus color-stealing.

In order to make a physical picture of $\mathfrak{P}_{\mathcal{F}}$ and thus
illustrate the tops function $\phi_{\mathcal{G}}\circ\tau_{\mathcal{F}}$ we
use a variant of the chaos game. To work at finite precision\ we partition the
set $\square\subset\mathbb{R}^{2}$ into a finite set of small rectangles, say
ten thousand of them, which we refer to as pixels. Each point $(x,y)\in
\square$ belongs to exactly one pixel, which we denote by $p((x,y))$.

Start from an arbitrary pair of points $(x_{0}^{\mathcal{F}},y_{0}%
^{\mathcal{F}})\in\square$ and $(x_{0}^{\mathcal{G}},y_{0}^{\mathcal{G}}%
)\in\square$. Let $K$ be a large number such as ten million. For $k=1,2,...K$
let $\sigma_{k}$ denote an element of $\{1,2,...,N\}$ chosen at random,
independently of all other choices. Let
\[
(x_{k}^{\mathcal{F}},y_{k}^{\mathcal{F}})=f_{\sigma_{k}}(x_{k-1}^{\mathcal{F}%
},y_{k-1}^{\mathcal{F}})\text{ and }(x_{k}^{\mathcal{G}},y_{k}^{\mathcal{G}%
})=g_{\sigma_{k}}(x_{k-1}^{\mathcal{G}},y_{k-1}^{\mathcal{G}})\text{.}%
\]
For each iterative step $k>100$, if the color of the pixel $p((x_{k}%
^{\mathcal{F}},y_{k}^{\mathcal{F}}))$ was not assigned at an earlier step
$l<k$ such that $\sigma_{l}\sigma_{l-1}\sigma_{l-2}...\sigma_{1}\overline
{1}>\sigma_{k}\sigma_{k-1}\sigma_{k-2}...\sigma_{1}\overline{1}$, then plot
the pixel $p((x_{k}^{\mathcal{F}},y_{k}^{\mathcal{F}}))$ in the color
$\mathfrak{P}_{\mathcal{G}}(x_{k}^{\mathcal{G}},y_{k}^{\mathcal{G}})$.

The reason this algorithm converges in practice to produce a stable physical
picture that approximates $\mathfrak{P}_{\mathcal{G}}\circ\phi_{\mathcal{G}%
}\circ\tau_{\mathcal{F}}$ is described in Chapter 4 of \cite{Bsuperfractals}.
Again, it depends on Birkhoff's ergodic theorem. Intuitively, ergodicity of
the shift transformation ensures that, almost always, the sequences
$\{(x_{k}^{\mathcal{F}},y_{k}^{\mathcal{F}})\}$ and $\{(x_{k}^{\mathcal{G}%
},y_{k}^{\mathcal{G}})\}$ repeatedly visit all of the pixels that represent
the points of $A_{\mathcal{F}}$ and $A_{\mathcal{G}}$ respectively. Let
$\sigma^{(k)}=\sigma_{k}\sigma_{k-1}\sigma_{k-2}...\sigma_{1}\overline{1}$.
Then the point $(x_{k}^{\mathcal{F}},y_{k}^{\mathcal{F}})$ is very close to
$\phi_{\mathcal{F}}(\sigma^{(k)})$ when $k$ is sufficiently large; indeed%
\[
d_{\mathbb{R}^{2}}(\phi_{\mathcal{F}}(\sigma^{(k)}),(x_{k}^{\mathcal{F}}%
,y_{k}^{\mathcal{F}}))\leq l^{k}d_{\mathbb{R}^{2}}(\phi_{\mathcal{F}%
}(\overline{1}),(x_{0}^{\mathcal{F}},y_{0}^{\mathcal{F}}))\text{.}%
\]
Similarly $(x_{k}^{\mathcal{G}},y_{k}^{\mathcal{G}})$ is very close to
$\phi_{\mathcal{G}}(\sigma^{(k)})$ when $k$ is sufficiently large. Hence, to a
good approximation, the color of the pixel $p(\phi_{\mathcal{F}}(\sigma
^{(k)}))$ is updated to become the color of the pixel $p(\phi_{\mathcal{G}%
}(\sigma^{(k)}))$ except when $\sigma^{(l)}>\sigma^{(k)}$ for some $l<k$ for
which $p(\phi_{\mathcal{F}}(\sigma^{(l)}))=p(\phi_{\mathcal{F}}(\sigma
^{(k)}))$. Let
\[
100<k_{1}<k_{2}<k_{3}<...<k_{M}\leq K
\]
denote the sequence of successive values of $k$ at which such updates occur.
Then $\{\sigma^{(k_{l})}\}_{l=1}^{M}$ is an increasing sequence of addresses,
each associated with a point in the pixel $p(\phi_{\mathcal{F}}(\sigma
^{(k_{1})}))$. Hence, again invoking ergodicity, $\{\sigma^{(k_{l})}%
\}_{l=1}^{M}$approaches the highest address of all points in the pixel
$p(\phi_{\mathcal{F}}(\sigma^{(k_{1})}))$. The address $\sigma^{(k_{M})}$ is
our approximation to $sup\{\tau_{\mathcal{F}}(\sigma):\sigma\in\tau
_{\mathcal{F}}(\phi_{\mathcal{F}}(\sigma^{(k_{1})}))\}$. In general we expect
it to become increasingly accurate with increasing $K$. According to this
approximation, the pixel $p(\phi_{\mathcal{F}}(\sigma^{(k_{1})}))$ is assigned
the colour of the pixel $p(\phi_{\mathcal{G}}(\tau_{\mathcal{F}}%
(\sigma^{(k_{M})})))$. Thus we obtain a sensible pixel-based approximation to
$\mathfrak{P}_{\mathcal{G}}\circ\phi_{\mathcal{G}}\circ\tau_{\mathcal{F}}$.

\begin{figure}
[ptb]
\begin{center}
\includegraphics[
height=1.9in,
width=4.8845in
]%
{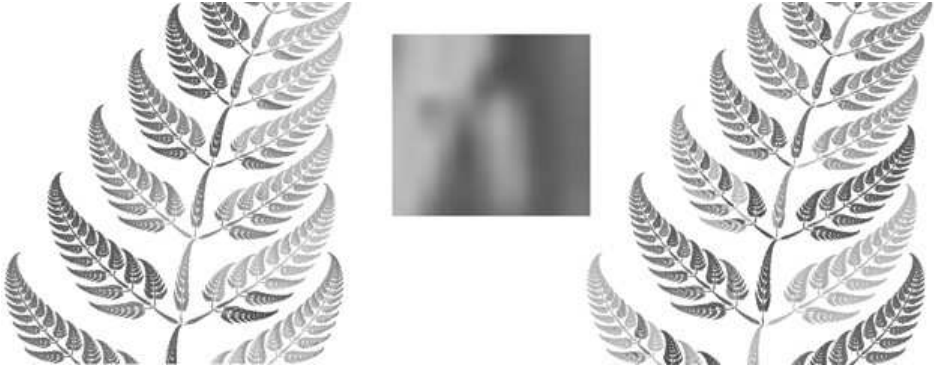}%
\caption{The ferns on the left and right are both obtained by fractal
transformations. The one on the left is continuous image of the central
image.}%
\label{ferns}%
\end{center}
\end{figure}

In Figure \ref{ferns} we illustrate two different fractal transformations from
a fern-like fractal to a filled square, computed using this algorithm. For the
picture on the left $\mathcal{F}$ and $\mathcal{G}$ are as discussed in
section \ref{example1}, with $\mathcal{C}_{\mathcal{F}}\prec\mathcal{C}%
_{\mathcal{G}}$, $\mathcal{C}_{\mathcal{F}}\neq\mathcal{C}_{\mathcal{G}}$, so
that $\phi_{\mathcal{G}}\circ\tau_{\mathcal{F}}$ is continuous. The picture
$\mathfrak{P}_{\mathcal{G}}$ is represented in the center of Figure
\ref{ferns}. It has been chosen to have apparently continuously varying
intensity so that the continuity of $\phi_{\mathcal{G}}\circ\tau_{\mathcal{F}%
}$ is illustrated by the smooth variation of intensity in the left-hand fern
image, which represents a close-up on $\mathfrak{P}_{\mathcal{F}}=$
$\mathfrak{P}_{\mathcal{G}}(\phi_{\mathcal{G}}\circ\tau_{\mathcal{F}})$. To
produce the picture on the right the IFS $\mathcal{G}$ has been switched, from
the one in Table \ref{CtsTable} to the one in Table \ref{DiscTable}, so that
$\phi_{\mathcal{G}}\circ\tau_{\mathcal{F}}$ is not continuous and
$\mathfrak{P}_{\mathcal{F}}(\phi_{\mathcal{G}}\circ\tau_{\mathcal{F}})$ is no
longer smoothly varying.%

\begin{figure}
[ptb]
\begin{center}
\includegraphics[
height=2.2987in,
width=3.0606in
]%
{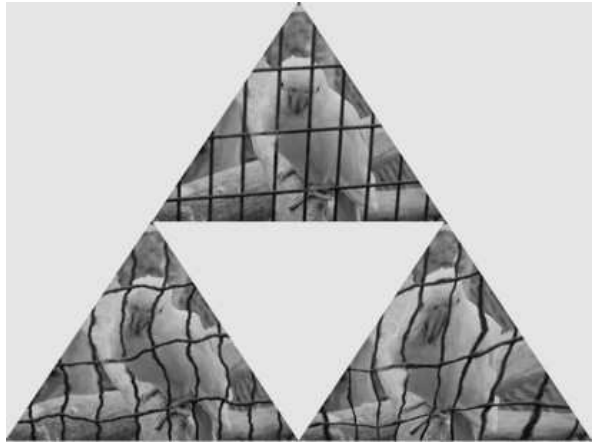}%
\caption{Two examples of fractal homeomorphisms applied to the picture at the
top. The tranformations from the top image to the one at bottom left is
area-preserving. }%
\label{threebirds}%
\end{center}
\end{figure}

In Figure \ref{threebirds} we illustrate two examples, computed using the
modified chaos game described here, in each of which the fractal
transformation $\phi_{\mathcal{G}}\circ\tau_{\mathcal{F}}:A_{\mathcal{F}%
}\rightarrow A_{\mathcal{G}}$ is a homeomorphism. The homeomorphisms are
constructed using IFSs of the form $\mathcal{F}_{\alpha,\beta,\gamma}$
discussed in sections \ref{hypersec} and \ref{example1}. In both examples
$A_{\mathcal{F}}=A_{\mathcal{G}}=\blacktriangle$, the filled triangle with
vertices at $A=(0,0),$ $B=(1,0)$, and $C=(0.5,1)$. Also in both cases,
$\mathfrak{P}_{\mathcal{G}}:\blacktriangle\rightarrow\mathfrak{C}$ corresponds
to the grayscale picture of a caged bird in the top triangle in Figure
\ref{threebirds}. The image at bottom left shows $\mathfrak{P}_{\mathcal{G}%
}\circ\phi_{\mathcal{G}}\circ\tau_{\mathcal{F}}$ when $\mathcal{F}%
=\mathcal{F}_{0.525,0.525,0.525}$ and $\mathcal{G}=\mathcal{F}%
_{0.475,0.475,0.475}$. In this case the corresponding subtriangles have the
same areas at all levels with the consequence that the fractal transformation
$\phi_{\mathcal{G}}\circ\tau_{\mathcal{F}}$ is area-preserving. To produce the
image at the bottom right we used $\mathcal{F}=\mathcal{F}_{0.4,0.6,0.475}$
and $\mathcal{G}=\mathcal{F}_{0.5,0.5,0.5}$.%

\section{\label{tdssec}The tops dynamical system}

In general, to determine the nature of the fractal transformation
$\phi_{\mathcal{G}}\circ\tau_{\mathcal{F}}:A_{\mathcal{F}}\rightarrow
A_{\mathcal{G}}$ we need to know the tops code space $\Omega_{\mathcal{F}}$.
Here we prove that $\Omega_{\mathcal{F}}$ is shift invariant. Consequently it
may be described in terms of the orbits of an associated dynamical system
$T_{\mathcal{F}}:A_{\mathcal{F}}\rightarrow A_{\mathcal{F}}$.

Throughout this section we assume that the transformations of the IFS
$\mathcal{F}$ are one-to-one. Let $\mathcal{S}_{\mathcal{F}}:\Omega
_{\mathcal{F}}\rightarrow\Omega$ denote the shift transformation, defined by
\[
\mathcal{S}_{\mathcal{F}}(\sigma_{1}\sigma_{2}\sigma_{3}...)=\sigma_{2}%
\sigma_{3}\sigma_{4}...
\]
for all $\sigma_{1}\sigma_{2}\sigma_{3}...\in\Omega_{\mathcal{F}}$. Let
\[
G_{\mathcal{F}}:=\{(x,\tau_{\mathcal{F}}(x)):x\in A_{\mathcal{F}}\}
\]
denote the graph of the the tops function $\tau_{\mathcal{F}}$.

\begin{lemma}
\label{lemma1}Let $(x,\sigma)\in G_{\mathcal{F}}$. Then $(f_{\sigma_{1}}%
^{-1}(x),\mathcal{S}_{\mathcal{F}}(\sigma))\in G_{\mathcal{F}}$.
\end{lemma}

\begin{proof}
$(x,\sigma)\in G_{\mathcal{F}}$ implies $x\in A_{\mathcal{F}}$, $\sigma
\in\Omega_{\mathcal{F}}$ and $\tau_{\mathcal{F}}(x)=\sigma$. In particular,
$\phi_{\mathcal{F}}(\sigma)=x$ for any $z\in\mathbb{X}$,
\[
\lim_{k\rightarrow\infty}f_{\sigma_{1}}\circ f_{\sigma_{2}}\circ
...f_{\sigma_{k}}(z)=x\text{.}%
\]
Using the continuity and invertibility of $f_{\sigma_{1}}$ it follows that
\[
\lim_{k\rightarrow\infty}f_{\sigma_{2}}\circ f_{\sigma_{3}}\circ
...f_{\sigma_{k}}(z)=f_{\sigma_{1}}^{-1}(x)\text{.}%
\]
This says that $\phi_{\mathcal{F}}(\mathcal{S}_{\mathcal{F}}(\sigma
))=f_{\sigma_{1}}^{-1}(x)$ which tells us that $\mathcal{S}_{\mathcal{F}%
}(\sigma)\in\phi_{\mathcal{F}}^{-1}(\{f_{\sigma_{1}}^{-1}(x)\})$.

Now suppose that there is $\omega\in\phi_{\mathcal{F}}^{-1}(\{f_{\sigma_{1}%
}^{-1}(x)\})$ with $\omega>\mathcal{S}_{\mathcal{F}}(\sigma)$. Then
$\phi_{\mathcal{F}}(\omega)=f_{\sigma_{1}}^{-1}(x)$ which implies
$f_{\sigma_{1}}(\phi_{\mathcal{F}}(\omega))=\phi_{\mathcal{F}}(\sigma
_{1}\omega)=x$. Let $\widetilde{\sigma}=\sigma_{1}\omega$. Then $\widetilde
{\sigma}>\sigma$ and $\phi_{\mathcal{F}}(\widetilde{\sigma})=x$ which
contradicts the assertion that $\sigma$ is the largest element of $\Omega$
such that $\phi_{\mathcal{F}}(\sigma)=x$. Hence $\mathcal{S}_{\mathcal{F}%
}(\sigma)\in\Omega_{\mathcal{F}}$ and $\tau_{\mathcal{F}}(f_{\sigma_{1}}%
^{-1}(x))=\mathcal{S}_{\mathcal{F}}(\sigma)$.
\end{proof}

\begin{lemma}
\label{lemma2}Let $(x,\sigma)\in G_{\mathcal{F}}$. Then $(f_{1}(x),1\sigma)\in
G_{\mathcal{F}}$.
\end{lemma}

\begin{proof}
$(x,\sigma)\in G_{\mathcal{F}}$ implies $\tau_{\mathcal{F}}(x)=\sigma$. Hence
$x=\phi_{\mathcal{F}}(\sigma)$ and so $f_{1}(x)=\phi_{\mathcal{F}}(1\sigma)$.

Now suppose that $(f_{1}(x),1\sigma)\notin G_{\mathcal{F}}$. Then there is
$\omega>1\sigma$ such that $\phi_{\mathcal{F}}(\omega)=f_{1}(x)$. But then
$\omega=1\widetilde{\sigma}$ where $\widetilde{\sigma}>\sigma$ and
$\phi_{\mathcal{F}}(1\widetilde{\sigma})=f_{1}(x)$. This implies
$\phi_{\mathcal{F}}(\widetilde{\sigma})=x$ with $\widetilde{\sigma}>\sigma$
which implies $\tau_{\mathcal{F}}(x)>\sigma$ which is a contradiction. Hence
$(f_{1}(x),1\sigma)\in G_{\mathcal{F}}$.
\end{proof}

It follows from Lemmas \ref{lemma1} and \ref{lemma2} that the mapping
$\widehat{T}_{\mathcal{F}}:G_{\mathcal{F}}\rightarrow G_{\mathcal{F}}$
specified by
\[
\widehat{T}_{\mathcal{F}}(x,\sigma)=(f_{\sigma_{1}}^{-1}(x),\mathcal{S}%
_{\mathcal{F}}(\sigma))\text{ for all }(x,\sigma)\in G_{\mathcal{F}}%
\]
is well-defined and onto.

In particular, the projection of $\widehat{T}_{\mathcal{F}}$ on $\Omega
_{\mathcal{F}}$ yields the symbolic dynamical system $\mathcal{S}%
_{\mathcal{F}}:\Omega_{\mathcal{F}}\rightarrow\Omega_{\mathcal{F}}$, because
from Lemma \ref{lemma2} we have
\[
\mathcal{S}_{\mathcal{F}}(\Omega_{\mathcal{F}})=\Omega_{\mathcal{F}}\text{.}%
\]
The projection of $\widehat{T}_{\mathcal{F}}:G_{\mathcal{F}}\rightarrow
G_{\mathcal{F}}$ onto $A_{\mathcal{F}}$ yields what we call the tops dynamical
system
\[
T_{\mathcal{F}}:A_{\mathcal{F}}\rightarrow A_{\mathcal{F}}%
\]
where
\begin{equation}
T_{\mathcal{F}}(x)=\left\{
\begin{array}
[c]{ccc}%
f_{1}^{-1}(x) & if & x\in D_{1}:=f_{1}(A_{\mathcal{F}}),\\
f_{2}^{-1}(x) & if & x\in D_{2}:=f_{2}(A_{\mathcal{F}})\backslash
f_{1}(A_{\mathcal{F}}),\\
. & . & .\\
f_{N}^{-1}(x) & if & x\in D_{N}:=f_{N}(A_{\mathcal{F}})\backslash%
{\textstyle\bigcup\limits_{n=1}^{N-1}}
f_{n}(A_{\mathcal{F}}),
\end{array}
\right.  \label{partition}%
\end{equation}
for all $x\in A_{\mathcal{F}}$. Lemma \ref{lemma2} implies
\[
T_{\mathcal{F}}(A_{\mathcal{F}})=A_{\mathcal{F}}\text{.}%
\]

\begin{theorem}
\label{TDSthm}The tops dynamical systems $T_{\mathcal{F}}:A_{\mathcal{F}%
}\rightarrow A_{\mathcal{F}}$ is related to the symbolic dynamical system
$\mathcal{S}_{\mathcal{F}}:\Omega_{\mathcal{F}}\rightarrow\Omega_{\mathcal{F}%
}$ by the tops function $\tau_{\mathcal{F}}:A_{\mathcal{F}}\rightarrow
\Omega_{\mathcal{F}}$, according to
\[
\mathcal{S}_{\mathcal{F}}=\tau_{\mathcal{F}}\circ T_{\mathcal{F}}\circ
\tau_{\mathcal{F}}^{-1}\text{.}%
\]
If $\Omega_{\mathcal{F}}\subset\Omega_{\mathcal{G}}$ then%
\[
(\phi_{\mathcal{G}}\circ\tau_{\mathcal{F}}\circ T_{\mathcal{F}}%
)(x)=(T_{\mathcal{G}}\circ\phi_{\mathcal{G}}\circ\tau_{\mathcal{F}})(x)\text{
for all }x\in A_{\mathcal{F}}\text{.}%
\]
If $\mathcal{C}_{\mathcal{F}}=\mathcal{C}_{\mathcal{G}}$ then the tops
dynamical systems $T_{\mathcal{F}}:A_{\mathcal{F}}\rightarrow A_{\mathcal{F}}$
and $T_{\mathcal{G}}:A_{\mathcal{G}}\rightarrow A_{\mathcal{G}}$ are
topologically conjugate.
\end{theorem}

\begin{proof}
Let $\Phi_{\mathcal{F}}=\tau_{\mathcal{F}}^{-1}$ be as discussed at the end of
section \ref{topsfnsec}. Then we claim that%
\[
T_{\mathcal{F}}\circ\Phi_{\mathcal{F}}=\Phi_{\mathcal{F}}\circ\mathcal{S}%
_{\mathcal{F}}\text{.}%
\]
Since $\mathcal{S}_{\mathcal{F}}$ maps $\Omega_{\mathcal{F}}$ onto itself and
$\Phi_{\mathcal{F}}$ maps $\Omega_{\mathcal{F}}$ onto $A_{\mathcal{F}}$ it
follows that the mapping $\Phi_{\mathcal{F}}\circ\mathcal{S}_{\mathcal{F}}$
takes $\Omega_{\mathcal{F}}$ onto $A_{\mathcal{F}}$. (Similarly,
$T_{\mathcal{F}}\circ\Phi_{\mathcal{F}}$ maps $\Omega_{\mathcal{F}}$ onto
$A_{\mathcal{F}}$.)

Let $\sigma=\sigma_{1}\sigma_{2}\sigma_{3}...\in\Omega_{\mathcal{F}}$. Then
$\mathcal{S}_{\mathcal{F}}(\sigma)=\sigma_{2}\sigma_{3}...\in\Omega
_{\mathcal{F}}$ and
\[
(\Phi_{\mathcal{F}}\circ\mathcal{S}_{\mathcal{F}})(\sigma)=\lim_{k\rightarrow
\infty}(f_{\sigma_{2}}\circ f_{\sigma_{3}}\circ...\circ f_{\sigma_{k}%
})(z)\text{.}%
\]
On the other hand
\begin{align*}
\Phi_{\mathcal{F}}(\sigma)  &  =\lim_{k\rightarrow\infty}(f_{\sigma_{1}}\circ
f_{\sigma_{2}}\circ...\circ f_{\sigma_{k}})(z)\\
&  =f_{\sigma_{1}}(\lim_{k\rightarrow\infty}(f_{\sigma_{2}}\circ f_{\sigma
_{3}}\circ...\circ f_{\sigma_{k}})(z))=f_{\sigma_{1}}(\Phi_{\mathcal{F}%
}(\sigma_{2}\sigma_{3}...))\text{,}%
\end{align*}
belongs to $A_{\mathcal{F}}$ and lies in the range of $f_{\sigma_{1}}$ and so
must belong to $D_{\sigma_{1}}$ as defined in Equation \ref{partition}. Hence
\[
(T_{\mathcal{F}}\circ\Phi_{\mathcal{F}})(\sigma)=\Phi_{\mathcal{F}}(\sigma
_{2}\sigma_{3}...)=(\Phi_{\mathcal{F}}\circ\mathcal{S}_{\mathcal{F}})(\sigma)
\]
for all $\sigma\in\Omega_{\mathcal{F}}$.

We now apply $\tau_{\mathcal{F}}$ to both sides of this last equation to
complete the proof of the first assertion in the theorem.

Now assume that $\Omega_{\mathcal{F}}\subset\Omega_{\mathcal{G}}$. Then,
since
\[
\mathcal{S}_{\mathcal{F}}(\sigma)=\mathcal{S}_{\mathcal{G}}(\sigma)
\]
for all $\sigma\in\Omega_{\mathcal{F}}$, it follows from the first part of the
theorem that%
\[
(\tau_{\mathcal{F}}\circ T_{\mathcal{F}}\circ\Phi_{\mathcal{F}})(\sigma
)=(\tau_{\mathcal{G}}\circ T_{\mathcal{G}}\circ\Phi_{\mathcal{G}})(\sigma)
\]
for all $\sigma\in\Omega_{\mathcal{F}}$. It follows that
\[
(\tau_{\mathcal{F}}\circ T_{\mathcal{F}}\circ\Phi_{\mathcal{F}}\circ
\tau_{\mathcal{F}})(x)=(\tau_{\mathcal{G}}\circ T_{\mathcal{G}}\circ
\Phi_{\mathcal{G}}\circ\tau_{\mathcal{F}})(x)
\]
for all $x\in A_{\mathcal{F}}$. But $\Phi_{\mathcal{F}}\circ\tau_{\mathcal{F}%
}=i_{A_{\mathcal{F}}}$ and
\[
(\Phi_{\mathcal{G}}\circ\tau_{\mathcal{F}})(x)=(\phi_{\mathcal{G}}\circ
\tau_{\mathcal{F}})(x)
\]
for all $x\in A_{\mathcal{F}}$. Hence%
\[
(\tau_{\mathcal{F}}\circ T_{\mathcal{F}})(x)=(\tau_{\mathcal{G}}\circ
T_{\mathcal{G}}\circ\phi_{\mathcal{G}}\circ\tau_{\mathcal{F}})(x)
\]
for all $x\in A_{\mathcal{F}}$. Applying $\phi_{\mathcal{G}}$ to both sides we
obtain%
\[
(\phi_{\mathcal{G}}\circ\tau_{\mathcal{F}}\circ T_{\mathcal{F}})(x)=(\phi
_{\mathcal{G}}\circ\tau_{\mathcal{G}}\circ T_{\mathcal{G}}\circ\phi
_{\mathcal{G}}\circ\tau_{\mathcal{F}})(x)\text{ }%
\]
for all $x\in A_{\mathcal{F}}$. But $\phi_{\mathcal{G}}\circ\tau_{\mathcal{G}%
}=i_{A_{\mathcal{G}}}$. This completes the proof of the second assertion in
the theorem.

Finally, let us suppose that $\mathcal{C}_{\mathcal{F}}=\mathcal{C}%
_{\mathcal{G}}$. Then Theorem 1 implies that $\phi_{\mathcal{G}}\circ
\tau_{\mathcal{F}}$ is a homeomorphism from $A_{\mathcal{F}}$ onto
$A_{\mathcal{G}}$. Also $\mathcal{C}_{\mathcal{F}}=\mathcal{C}_{\mathcal{G}}$
implies $\Omega_{\mathcal{F}}=\Omega_{\mathcal{G}}$ which implies, via the
previously proven part of this theorem,
\[
T_{\mathcal{F}}(x)=(\phi_{\mathcal{G}}\circ\tau_{\mathcal{F}})^{-1}\circ
T_{\mathcal{G}}\circ\phi_{\mathcal{G}}\circ\tau_{\mathcal{F}})(x)
\]
for all $x\in A_{\mathcal{F}}$.
\end{proof}

If the domains $\{D_{n}:n=1,2,..,N\}$ are known then it is easy to compute the
tops function. Just follow the orbit of $x$ under the tops dynamical system
and keep track of the sequence of indices $\sigma_{1}\sigma_{2}\sigma
_{3}\sigma_{4}...$ \ visited by the orbit.

In the special case where the IFS is totally disconnected and the $f_{n}$s are
one-to-one then $T:A\rightarrow A$ is defined by $T(x)=f_{n}^{-1}(x)$ where
$n$ is the unique index such that $x\in f_{n}(A)$. This dynamical system has
been considered elsewhere, for example in \cite{BaDe} and \cite{kieninger}. In
this case $\phi_{\mathcal{F}}:\Omega\rightarrow A_{\mathcal{F}}$ is a
homeomorphism, $\tau_{\mathcal{F}}=\phi_{\mathcal{F}}^{-1}$, and
$T_{\mathcal{F}}$ it is conjugate to the shift transformation according
to$T_{\mathcal{F}}=\phi_{\mathcal{F}}\circ\mathcal{S}_{\mathcal{F}}\circ
\phi_{\mathcal{F}}^{-1}$.

Theorem \ref{TDSthm} says\ in particular that $T_{\mathcal{F}}:A_{\mathcal{F}%
}\rightarrow A_{\mathcal{F}}$ is a factor of $\mathcal{S}_{\mathcal{F}}%
:\Omega_{\mathcal{F}}\rightarrow\Omega_{\mathcal{F}}$, and as defined for
example in \cite{katok} p.68, because $\Phi_{\mathcal{F}}\circ\mathcal{S}%
_{\mathcal{F}}=T_{\mathcal{F}}\circ\Phi_{\mathcal{F}}$ where $\Phi
_{\mathcal{F}}=\tau_{\mathcal{F}}^{-1}$ is continuous; this tells us that the
topological entropy of $T_{\mathcal{F}}$ is less than or equal to the
topological entropy of $\mathcal{S}_{\mathcal{F}}$, \cite{katok} Proposition
3.1.6, p.111. If $\tau_{\mathcal{F}}$ is continuous then Theorem \ref{TDSthm}
says that the two dynamical systems $T_{\mathcal{F}}:A_{\mathcal{F}%
}\rightarrow A_{\mathcal{F}}$ and $\mathcal{S}_{\mathcal{F}}:\Omega
_{\mathcal{F}}\rightarrow\Omega_{\mathcal{F}}$ are topologically conjugate,
see \cite{katok} p. 60, and it follows that the two systems must have the same
topological entropy.

This suggests that we may compare the complexity of some subsets of
$\mathbb{R}^{2}$ by assigning to them the topological entropy of a
corresponding shift dynamical system. Let $\mathcal{M}$ denote the set of all
attractors of hyperbolic IFSs in $\mathbb{R}^{2}$, whose transformations are
all affine and invertible, such that the associated tops function is
continuous. Then we can define the topological entropy of each $A_{\mathcal{F}%
}$ to be the infimum of the entropies of the set of corresponding shift
dynamical systems. In this way we arrive at a geometry-based definition of the
topological entropy of some subsets of $\mathbb{R}^{2}$. Is it useful?

\end{document}